\documentclass{article}
\usepackage{amsmath}
\usepackage{amsfonts}

\newtheorem{theorem}{Theorem}

\newtheorem{lemma}[theorem]{Lemma}

\input{tcilatex}
\begin{document}

\title{A note on a global invertibility of mappings on $R^{n}$}
\author{ Marek Galewski }
\maketitle
\date{}

\begin{abstract}
\noindent \noindent We provide sufficient conditions for a mapping $%
f:R^{n}\rightarrow R^{n}$ to be a global diffeomorphism in case it is
strictly (Hadamard) differentiable. We use classical local invertibility
conditions together with the non-smooth critical point theory.

\noindent \textbf{Math Subject Classifications}: 57R50, 58E05

\end{abstract}

\section{Introduction}

In this note we consider locally invertible mappings $f:R^{n}\rightarrow
R^{n}$ that are strictly differentiable and which need not be assumed
continuously Fr\'{e}chet - differentiable. We are interested in sufficient
conditions for $f$ to be a global diffeomorphism. The notion of strict
differentiability (Hadamard derivative) is intermediate between Fr\'{e}chet
- differentiability and $f$ being a $C^{1}$ mapping and it reads as follows,
see \cite{ochalbook}, $f:D\rightarrow R^{n}$ defined on $D\subset R^{n}$ is
Hadamard (strictly) differentiable at $x_{0}\in D$, if there exists $%
f^{\prime }(x_{0})\in R^{n}$ such that%
\begin{equation*}
\lim_{w\rightarrow x_{0},t\rightarrow 0^{+}}\frac{f(w+tz)-f(w)}{t}%
=\left\langle f^{\prime }(x_{0}),z\right\rangle \text{ for all }z\in R^{n}
\end{equation*}%
provided the convergence is uniform for $z$ in compact sets. \medskip

Using non-smooth critical point theory applied to a functional $x\rightarrow 
\frac{1}{2}\left\Vert f\left( x\right) \right\Vert ^{2}$, where $\left\Vert
\cdot \right\Vert $ stands for the Euclidean norm, we provide sufficient
conditions for $f$ to be global diffeomorphism. Since $f$ is not $C^{1},$
functional $x\rightarrow \frac{1}{2}\left\Vert f\left( x\right) \right\Vert
^{2}$ need not be such and thus we cannot use a smooth result proved via a
classical mountain pass theorem by Katriel \cite{katriel}

\begin{theorem}
\label{MainTheo copy(2)}Let $X,$ $B$ be finite dimensional Euclidean spaces.
Assume that $f:X\rightarrow B$ is a $C^{1}$-mapping such that\newline
(a1) $f^{\prime }(x)$ is invertible for any $x\in X$;\newline
(a2) $\left\Vert f\left( x\right) \right\Vert \rightarrow \infty $ as $%
\left\Vert x\right\Vert \rightarrow \infty $\newline
then $f$ is a diffeomorphism.
\end{theorem}

Since the mountain geometry in the non-smooth setting will work in this
case, we shall use the ideas of Katriel in our reasoning together with some
ideas from \cite{SIW}, where an infite dimensional version of the above
result is to be found. The local invertibility results we base on are as
follows, \cite{rad}, see also \cite{guttierez}.

\begin{lemma}
\label{lem1}Let $D$ will be an open set of $R^{n}$ and let $f:D\rightarrow
R^{n}$ be a Fr\'{e}chet - differentiable map and the \ following condition
holds: $\det [f^{\prime }(x)]\neq 0$ for every $x\in D.$Then $f$ is a local
diffeomorphism.
\end{lemma}

For background on non-smooth analysis we refer to \cite{CLARKE2} and in \cite%
{ochalbook}. If $f:R^{n}\rightarrow R^{n}$ is a locally Lipschitz continuous
function for $u,z\in R^{n}$ define $f^{0}(u;z)$ for the generalized
directional derivative at the point $u$ along the direction $z$ by 
\begin{equation*}
f^{0}(u;z):=\limsup_{w\rightarrow u,\,t\rightarrow 0^{+}}\frac{f(w+tz)-f(w)}{%
t}\,.
\end{equation*}%
The generalized gradient of $J$ in $u$ is the set 
\begin{equation*}
\partial f(u):=\{\xi \in R^{n}:\langle \xi ,z\rangle \leq f^{0}(u;z),\;\text{%
for all }\,z\in R^{n}\}.
\end{equation*}

If $f$ is strictly differentiable the Clarke subdifferential reduces to a
singleton, i.e. its derivative.

A point $u$ is called a critical point of $f,$ if $0\in \partial f(u)$. A
locally Lipschitz continuous functional $J:R^{n}\rightarrow \mathbb{R}$ is
said to fulfill the non-smooth Palais-Smale condition if every sequence $%
\{u_{n}\}$ in $R^{n}$ such that $\{J(u_{n})\}$ is bounded and 
\begin{equation*}
J^{0}(u_{n};u-u_{n})\geq -\varepsilon _{n}\Vert u-u_{n}\Vert
\end{equation*}%
\ for all $u\in R^{n}$, where $\varepsilon _{n}\rightarrow 0^{+},$ admits a
convergent subsequence.

We will use the following version of the mountain pass lemma in the
nonsmooth setting.

\begin{theorem}
\label{MPT} \cite{MoVa}Let $J:R^{n}\rightarrow \mathbb{R}$ be a locally
Lipschitz continuous functional satisfying the non-smooth Palais-Smale
condition. If there exist $u_{1},u_{2}\in R^{n}$, $u_{1}\neq u_{2}$ and $%
r\in (0,\left\Vert u_{2}-u_{1}\right\Vert )$ such that 
\begin{equation*}
\inf \{J(u):\left\Vert u-u_{1}\right\Vert =r\}\geq \max \{J(u_{1}),J(u_{2})\}
\end{equation*}%
and we denote by $\Gamma $ the family of continuous paths $\gamma
:[0,1]\rightarrow R^{n}$ joining $u_{1}$ and $u_{2},$ then%
\begin{equation*}
c:=\underset{\gamma \in \Gamma }{\inf }\underset{s\in \lbrack 0,1]}{\max }%
J(\gamma (s))\geq \max \{J(u_{1}),J(u_{2})\}
\end{equation*}%
is a critical value for $R^{n}$\ and $K_{c}\backslash \{u_{1},u_{2}\}\neq
\emptyset $, where $K_{c}$ is the set of critical points at the level $c$.
\end{theorem}

\section{Results}

\begin{theorem}
\label{MainTheo} If $f:R^{n}\rightarrow R^{n}$ is a strictly differentiable
mapping such that\newline
(b1) for any $y\in R^{n}$ the functional $\varphi :R^{n}\rightarrow 
%TCIMACRO{\U{211d} }%
%BeginExpansion
\mathbb{R}
%EndExpansion
$ defined by 
\begin{equation*}
\varphi \left( x\right) =\frac{1}{2}\left\Vert f\left( x\right)
-y\right\Vert ^{2}
\end{equation*}%
is coercive, i.e. $\varphi \left( x\right) \rightarrow \infty $ as $%
\left\Vert x\right\Vert \rightarrow \infty $ for any $y\in R^{n}$; \newline
(b2) for any $x\in R^{n}$ we have $\det [f^{\prime }(x)]\neq 0$ \newline
then $f$ is a diffeomorphism.
\end{theorem}

By Lemma \ref{lem1} condition (\textit{b2}) implies that $f$ defines a local
diffeomorphism. Thus it is sufficient to show that $f$ is onto and one to
one.\bigskip

Let us fix any point $y\in R^{n}$. Observe that $\varphi $ being a
composition of a $C^{1}$ mapping and a strictly differentiable mapping is
strictly differentiable and therefore it is locally Lipschitz continuous
with $\partial f(x)$ being equal to $(\left\{ f\left( x\right) -y)\circ
f^{\prime }(x)\right\} $ for any $x\in R^{n}$. Since $\varphi $ is
continuous and coercive it has an argument of a minimum $\overline{x}$,
which satisfies the non-smooth Fermat's rule, i.e. 
\begin{equation*}
0\in \left\{ f\left( \overline{x}\right) -y)\circ f^{\prime }(\overline{x}%
)\right\}
\end{equation*}%
which means that $0=f^{\prime }(\overline{x})^{T}f\left( \overline{x}\right)
-y)$, where $\xi ^{T}$ denotes the transpose of the matrix $\xi $. Since by
(b2) $\det f^{\prime }(\overline{x})\neq 0$, we see that $f\left( \overline{x%
}\right) -y=0$. Thus $f$ is surjective.\bigskip

Now we argue by contradiction that $f$ is one to one.\textbf{\ }Suppose
there are $x_{1}$ and $x_{2}$, $x_{1}\neq x_{2}$, $x_{1}$, $x_{2}\in R^{n}$,
such that $f\left( x_{1}\right) =f\left( x_{2}\right) =a\in R^{n}$. We will
apply Theorem \ref{MPT}. We put $e=x_{1}-x_{2}$ and define mapping $%
g:R^{n}\rightarrow R^{n}$ and a locally Lipschitz functional $\psi
:R^{n}\rightarrow 
%TCIMACRO{\U{211d} }%
%BeginExpansion
\mathbb{R}
%EndExpansion
$ by 
\begin{equation*}
g\left( x\right) =f\left( x+x_{2}\right) -a\text{ and }\psi \left( x\right) =%
\frac{1}{2}\left\Vert g\left( x\right) \right\Vert ^{2}.
\end{equation*}%
Note that $\psi \left( e\right) =\psi \left( 0\right) =0.$ By (\textit{b1}) $%
\varphi $ is coercive, so it satisfies the non-smooth\ Palais-Smale
condition. The same conclusion holds for functional $\psi $. Fix $\rho >0$
such that $\rho <\left\Vert e\right\Vert .$ By the classical Weierstrass
Theorem $\psi $ has an argument of a minimum over $\partial \overline{%
B\left( 0,\rho \right) }$ which we denote by $w$ and which is non-zero and
different from $x_{1}$ and $x_{2}$. Thus $\psi \left( w\right) >0$.
Therefore 
\begin{equation*}
\inf_{\left\Vert x\right\Vert =r}\psi (x)\geq \psi \left( w\right) >0=\psi
\left( e\right) =\psi \left( 0\right) .
\end{equation*}%
Thus by Theorem \ref{MPT} applied to $J=\psi $ we note that $\psi $ has a
critical point $v\neq 0$, $v\neq e$ and such that 
\begin{equation*}
\psi ^{^{\prime }}(v)=(f\left( v+x_{2}\right) -a)\circ f^{\prime
}(v+x_{2})=0.
\end{equation*}%
Since $\det f^{\prime }(v+x_{2})\neq 0$ we see that $f\left( v+x_{2}\right)
-a=0$. This means that either $v=0$ or $v=e$. Thus we obtain a contradiction
which shows that $f$ is a one to one operator.

It remains to comment on a situation when it is convenient to equip $R^{n}$
with a norm making it a Banach space without Euclidean structure, for
example the $p-$norm. These ideas come from \cite{GGS}.

\begin{theorem}
\label{MainTheo copy(1)}Assume that $f:R^{n}\rightarrow R^{n}$ is a strictly
differentiable mapping, $\eta :R^{n}\rightarrow 
%TCIMACRO{\U{211d} }%
%BeginExpansion
\mathbb{R}
%EndExpansion
_{+}$ is a $C^{1}$ functional and that the following conditions hold\newline
(c1) $\left( \eta \left( x\right) =0\Longleftrightarrow x=0\right) $ and $%
\left( \eta ^{^{\prime }}\left( x\right) =0\Longleftrightarrow x=0\right) $; 
\newline
(c2) for any $y\in R^{n}$ the functional $\varphi :X\rightarrow 
%TCIMACRO{\U{211d} }%
%BeginExpansion
\mathbb{R}
%EndExpansion
$ given by the formula 
\begin{equation*}
\varphi \left( x\right) =\eta \left( f\left( x\right) -y\right)
\end{equation*}%
is coercive;\newline
(c3) $\det f^{\prime }(x)\neq 0$ for any $x\in R^{n}$,\newline
(c4) there exist positive constants $\alpha $, $\beta $, $M$ such that 
\begin{equation*}
\eta \left( x\right) \geq c\left\Vert x\right\Vert ^{\alpha }\text{ for }%
\left\Vert x\right\Vert \leq M.
\end{equation*}%
\newline
Then $f$ is a diffeomorphism from $R^{n}$ onto $R^{n}$.
\end{theorem}

To indicate some differences in the proof note that by Fermat's Principle we
get 
\begin{equation*}
0\in \left\{ \eta ^{^{\prime }}(f\left( \overline{x}\right) -y)\circ
f^{\prime }(\overline{x})\right\} \text{.}
\end{equation*}%
By (\textit{c3)} and (c1) we get $f\left( \overline{x}\right) -y=0$. Define $%
g$ as before and put $\psi \left( x\right) =\eta \left( g\left( x\right)
\right) $. By $w$ we denote again an argument of a minimum of $\psi $ over $%
\partial \overline{B\left( 0,\rho \right) }$ with$\rho <\min \left\{
M,\left\Vert e\right\Vert \right\} $. By (c4) we get 
\begin{equation*}
\inf_{\left\Vert x\right\Vert =\rho }\psi (x)>0=\psi \left( e\right) =\psi
\left( 0\right) .
\end{equation*}%
The remaining part of the proof follows in the same manner.

\begin{tabular}{l}
Marek Galewski \\ 
Institute of Mathematics, \\ 
Lodz University of Technology, \\ 
Wolczanska 215, 90-924 Lodz, Poland, \\ 
marek.galewski@p.lodz.pl%
\end{tabular}

\end{document}